\documentclass{amsart}
\usepackage{amsfonts}
\newtheorem{theorem}{Theorem}

\newtheorem{lm}{Lemma}

\newtheorem{cor}{Corollary}
\newtheorem{prop}{Proposition}
\newtheorem{nota}{Notation}

\begin{document}
\title{A refined realization theorem in the context of the Schur-Szeg\H{o} 
composition}
\author{Vladimir P. Kostov}
\address{Universit\'e de Nice, Laboratoire de 
Math\'ematiques, Parc Valrose, 06108 Nice, Cedex 2, France}
\email{kostov@math.unice.fr}

\thanks{Research supported by the French Foundation CNRS under Project 20682.}

\subjclass[2000]{Primary 12A10 Secondary 30D99}

\keywords{Schur-Szeg\H{o} composition; composition factor; entire function}

\maketitle

\begin{abstract}
Every polynomial of the form $P=(x+1)(x^{n-1}+c_1x^{n-2}+\cdots +c_{n-1})$ is representable as Schur-Szeg\H{o} composition of $n-1$ polynomials of the form $(x+1)^{n-1}(x+a_i)$, where the numbers $a_i$ are unique up to permutation. We give necessary and sufficient conditions upon the possible values of the $8$-vector whose components are the number of positive, zero, negative and complex roots of a real polynomial $P$ and the number of positive, zero, negative and complex among the quantities $a_i$ corresponding to $P$. A similar result is proved about entire functions of the form $e^xR$, where $R$ is a polynomial.

\end{abstract}

\section{Introduction}
\subsection{Schur-Szeg\H{o} composition and the mapping $\Phi$}

In the present paper we prove a realization theorem in the context of the {\em Schur-Szeg\H{o} composition} of polynomials. For the two polynomials of degree $n$, $A:=\sum _{j=0}^n{n\choose j}a_jx^j$ and $B:=\sum _{j=0}^n{n\choose j}b_jx^j$, their Schur-Szeg\H{o} composition is defined by the formula 

$$A*B:=\sum _{j=0}^n{n\choose j}a_jb_jx^j~.$$ 
This formula is valid for any complex polynomials. In this paper we are interested mainly in the case when $A$ and $B$ are real.

Observe that when one considers the polynomials as degree $n+k$ ones, with $k$ leading zero coefficients, then the formula for their Schur-Szeg\H{o} composition will be a different one. To avoid such ambiguity, we assume throughout the paper that the leading coefficient of at least one of the composed polynomials is nonzero. 

The polynomial $(x+1)^n$ plays the role of unity in the sense that 

\begin{equation}\label{identity}
(x+1)^n*A=A
\end{equation}
for any polynomial $A$. Schur-Szeg\H{o} composition of polynomials is commutative and associative. It can be defined for an arbitrary number of polynomials by the formula

$$ A_1*\cdots *A_s:=\sum _{j=0}^n{n\choose j}a_j^1\cdots a_j^sx^j~,~~{\rm where}~~A_i:=\sum _{j=0}^n{n\choose j}a_j^ix^j~.$$
The reader can find a more detailed information about the Schur-Szeg\H{o} composition in the monographies \cite{Pr} and \cite{RS}. 
\newpage

The Schur-Szeg\H{o} composition gives rise to a mapping $\Phi$ (defined below) from the space of polynomials of degree $n-1$ into itself in the following way. We consider polynomials of degree $n$ having one of their roots at $(-1)$:  

\begin{equation}\label{Pn-1}
 P:=(x+1)(x^{n-1}+c_1x^{n-2}+\cdots +c_{n-1})
\end{equation}
Each such polynomial is representable as Schur-Szeg\H{o} composition of $n-1$ polynomials (called {\em composition factors}) of the form $K_{a_i}:=(x+1)^{n-1}(x+a_i)$, where the complex numbers $a_i$ are uniquely defined. In the case of real polynomials part 
of the numbers $a_i$ are real, while the rest form complex conjugate couples. This result has been announced in \cite{Ko1} and proved in \cite{AlKo}. 

Notice that the numbers $(-a_i)$ can be viewed as roots of another polynomial. In this sense we obtain a mapping from the space of degree $n-1$ polynomials into itself in the following way. Denote by $\sigma _j$ the elementary symmetric polynomials of the quantities $a_i$, i.e. $\sigma _j:=\sum _{1\leq i_1<i_2<\cdots <i_j\leq n-1}a_{i_1}\cdots a_{i_j}$. The mapping $\Phi$ is defined like this:

$$\Phi ~: ~(c_1,\ldots ,c_{n-1})\mapsto (\sigma _1,\ldots ,\sigma _{n-1})~.$$
The mapping $\Phi$ is affine non-degenerate (see \cite{AlKo}) and its eigenvalues are rational positive numbers (see \cite{Ko2}). For other properties of this mapping see \cite{Ko4} and \cite{KoShMa}. Denote in the case of a real polynomial $P$ by $\rho$ the number of the real roots of the polynomial $P/(x+1)$ and by $r$ the number of the real among the quantities $a_i$. Notice that $[(n-1-\rho )/2]$ is the number of complex conjugate couples of roots of $P/(x+1)$ ($[.]$ stands for the integer part of) and $[(n-1-r)/2]$ is the number of complex conjugate couples of quantities $a_j$. A priori, $0\leq \rho ,r\leq n$ and the parity of the numbers $r$ and $\rho$ must be the same. In paper \cite{Ko3} the following question is asked: 

{\em When these natural restrictions are respected, what can be the values of the couple $(r,\rho )$?} 

The answer given there is: 

{\em All possible values are attained at some polynomials having all their roots distinct and for which the corresponding quantities $a_i$ are also distinct.} 

In other words, all a priori admissible couples $(r,\rho )$ are realizable. (This is a realization theorem.)

A similar realization theorem has been proved in \cite{Ko6} for an analog of the Schur-Szeg\H{o} composition in the case of entire transcendental functions. Consider the two transcendental functions $f$ and $g$ represented by the convergent everywhere in $\mathbb{C}$ series $\sum _{j=0}^{\infty}\gamma _jx^j/j!$ and $\sum _{j=0}^{\infty}\delta _jx^j/j!$, respectively. Their Schur-Szeg\H{o} composition is defined by the formula 

$$f*g:=\sum _{j=0}^{\infty}\gamma _j\delta _jx^j/j!~.$$ 
As in the case of polynomials, Schur-Szeg\H{o} composition is commutative, associative and can be defined for any finite number of entire transcendental functions by the formula

$$ f_1*\cdots *f_s:=\sum _{j=0}^{\infty}\gamma _j^1\cdots \gamma _j^sx^j/j!~,~~{\rm where}~~f_i:=\sum _{j=0}^{\infty}\gamma _j^ix^j/j!~.$$

Consider transcendental functions of the form $e^xR$, where $R$ is a polynomial of degree $n-1$, $R(0)=1$. Such a function is representable as Schur-Szeg\H{o} composition 

\begin{equation}\label{compex}
 e^xR=\kappa _{a_1}*\cdots *\kappa _{a_{n-1}}~,
\end{equation}
where the {\em composition factors} $\kappa _{a_i}$ are of the form 
$e^x(1+x/a_j)$, the numbers $a_j$ being uniquely defined up to permutation. To extend the formula to the case when $P(0)$ is not necessarily $1$, one has to admit the presence of composition factors $e^xc$, $c\neq 0$, and $e^xx$ (one can say that the latter correspond to the case $a_i=0$). When the polynomial $P$ is real, part of the numbers $a_j$ are real and the rest form complex conjugate couples. 

Set $\tilde{\sigma}_j:=\sum _{1\leq i_1<i_2<\cdots <i_j\leq n-1}1/a_{i_1}\cdots a_{i_j}$ and 
$R:=\tilde{c}_{n-1}x^{n-1}+\tilde{c}_{n-2}x^{n-2}+\cdots +\tilde{c}_1x+1$. For transcendental functions the mapping $\Phi$ is defined by the formula

$$ \Phi ~:~(\tilde{c}_1,\ldots ,\tilde{c}_{n-1})\mapsto (\tilde{\sigma}_1,\ldots ,\tilde{\sigma}_{n-1})~.$$ 
As in the case of polynomials this mapping is affine. (This follows from the formulae 
$e^x*f(x)=f(x)$ and $xe^x*f(x)=xf'(x)$ 
which are readily checked; $f$ denotes an entire function. The formulae imply that every coefficient of $R$ is a polynomial in the quantities $1/a_i$. This polynomial is symmetric and in each monomial each factor $1/a_i$ appears in degree $0$ or $1$.) It is shown in \cite{Ko6} that all a priori admissible couples $(r,\rho )$ are realizable, where the quantities $r$ and $\rho$ are defined by analogy with the case of polynomials.

\subsection{The new results}

When Schur-Szeg\H{o} composition of real polynomials is considered, then it is important to distinguish the real positive, negative and zero roots of the composed polynomials. (About the number of real positive or negative roots of the Schur-Szeg\H{o} composition of two hyperbolic or real polynomials see respectively \cite{Ko1} and \cite{Ko5}.) Propositions~\ref{necessarypol} and \ref{necessaryex} show that the properties of the mapping $\Phi$ (defined for polynomials or for entire functions) are not the same with regard to the real positive and the real negative roots. We use the following notation:

\begin{nota}
{\rm We set $b_j:=-j/(n-j)$ for $j=0,1,\ldots ,n-1$, $b_n=-\infty$.}
\end{nota}

\begin{prop}\label{necessarypol}


{\rm (1)} If the polynomial $P$ from (\ref{Pn-1}) has $m$ positive roots counted with multiplicity ($m\geq 0$) and a $k$-fold root at $0$ ($k\geq 0$), then there are at least $m+\max (0,k-1)$ negative and distinct among the numbers $a_i$ out of which $\max (0,k-1)$ equal $b_1$, $\ldots$, $b_{k-1}$; if  $k\geq 1$, then one of the numbers $a_i$ equals $0$.

{\rm (2)} If there are $q$ numbers $a_i$ which equal $0$ and $q_1$ which are positive, then the polynomial $P$ has at least $q_1+\max (0,q-1)$ negative roots counted with multiplicity; for $q\geq 1$ it has a root at $0$.
\end{prop}

\begin{prop}\label{necessaryex}


{\rm (1)} If the polynomial $R$ has $m$ positive roots counted with multiplicity and a $k$-fold root at $0$, then there are at least $m+\max (0,k-1)$ negative and distinct among the numbers $a_i$ out of which $\max (0,k-1)$ equal $(-j)$, $j=1,\ldots ,k-1$. For $k\geq 1$ one composition factor equals $e^xx$.

{\rm (2)} If there are $q$ numbers $a_i$ which equal $0$ and $q_1$ which are positive, then the polynomial $R$ has at least $q_1+\max (0,q-1)$ negative roots counted with multiplicity; for $q\geq 1$ it has a root at $0$.
\end{prop}

The propositions are proved in Section~\ref{proofsprops}. In their proofs we use some facts about the Schur-Szeg\H{o} composition, see Section~\ref{knownfacts}. 

Theorems~\ref{maintm1} and \ref{maintm2} below show that the necessary conditions expressed by the propositions are in fact sufficient as well. In this sense they are realization theorems. Before formulating the theorems we analyse in detail the possible number of real positive, negative and zero among the roots of the polynomial $P$ or $R$ on the one-hand side and the numbers $(-a_i)$ on the other. We use the same notation as in the propositions.

Suppose that $k\geq 1$ (hence $q\geq 1$) and that among the numbers $(-a_i)$ corresponding to the polynomial $P$ or $R$ there are $q_1$ negative ones and $q_{\mathbb{C}}/2$ complex conjugate couples. Hence the quantity of positive numbers $(-a_i)$ is 

$$k-1+m+r,{\rm ~for~some~}r\geq 0~{\rm (see~part~(1)~of~the~propositions).}$$ 
There are $q-1+q_1+s$ negative among the roots of the polynomial (where $s\geq 0$, see part (2) of the propositions) and $k_{\mathbb{C}}/2$ complex conjugate couples. If $k=0$, then $q=0$ and the above two numbers equal $m+r$ and $q_1+s$.

We distinguish four cases. In Cases 1) and 2) (respectively 3) and 4)) we suppose that $k_{\mathbb{C}}\leq r$ (respectively $k_{\mathbb{C}}>r$). 
One has $k=1$ in Cases 1) and 3) and $k=0$ in Cases 2) and 4).

{\em Case 1)} Set $r=k_{\mathbb{C}}+k_1$. Hence $s=q_{\mathbb{C}}+k_1$. We present the situation schematically like this:

\[ \begin{array}{rcccccc}
&&&&&k_{\mathbb{C}}&\\ 
{\rm roots~of~}P/(x+1)~{\rm or~}R&&q-1+q_1+\overbrace{q_{\mathbb{C}}+k_1}^s&&k&&m\\ 
&&- - - - - - -&&0&&+ + + + + + +\\ 
{\rm numbers}~(-a_i)&&q_1&&q&&k-1+m+\underbrace{k_{\mathbb{C}}+k_1}_r\\&&&&&q_{\mathbb{C}}&\end{array}
\] 
The sequences of minus and plus signs to the left and right in the middle line symbolize the negative and positive half-axes; the numbers $k_{\mathbb{C}}$ and $q_{\mathbb{C}}$ are put away from that line, i.e. ``away from the real axis'' 
because these are the quantities of complex (not real) roots.

{\em Case 2)} is defined by the conditions $k=0$ (hence $q=0$), $k_{\mathbb{C}}\leq r$. The case can be presented like this:

$$ \begin{array}{rcccccc}
&&&&&k_{\mathbb{C}}&\\ 
{\rm roots~of~}P/(x+1)~{\rm or~}R&&q_1+\overbrace{q_{\mathbb{C}}+k_1}^s&&0&&m\\ 
&&- - - - - - -&&0&&+ + + + + + +\\ 
{\rm numbers}~(-a_i)&&q_1&&0&&m+\underbrace{k_{\mathbb{C}}+k_1}_r\\&&&&&q_{\mathbb{C}}& \end{array}
$$

In {\em Case 3)}, i.e. when $k\geq 1$ (hence $q\geq 1$) and $k_{\mathbb{C}}>r$, one can set $k_{\mathbb{C}}:=r+\delta$, $q_{\mathbb{C}}:=s+\delta$ and the situation admits the following presentation: 

$$ \begin{array}{rcccccc}
&&&&&r+\delta&\\ 
{\rm roots~of~}P/(x+1)~{\rm or~}R&&q-1+q_1+s&&k&&m\\ 
&&- - - - - - -&&0&&+ + + + + + +\\ 
{\rm numbers}~(-a_i)&&q_1&&q&&k-1+m+r\\&&&&&s+\delta&\end{array}
$$

Finally, in {\em Case 4)}, i.e. when $k=0=q$ and $k_{\mathbb{C}}>r$, the presentation looks like this:
$$ \begin{array}{rcccccc}
&&&&&r+\delta&\\ 
{\rm roots~of~}P/(x+1)~{\rm or~}R&&q_1+s&&0&&m\\ 
&&- - - - - - -&&0&&+ + + + + + +\\ 
{\rm numbers}~(-a_i)&&q_1&&0&&m+r\\&&&&&s+\delta&\end{array}
$$

We say that a polynomial $P$ or $R$ realizes Case 1), 2), 3) or 4) if the number of its positive, zero, negative and complex roots and these numbers defined for the quantities $(-a_i)$ are as shown on the above figures, and if the non-zero roots are distinct and the non-zero quantities $a_i$ are also distinct.

\begin{theorem}\label{maintm1}
Cases 1) -- 4) are realizable by some polynomials $P$.
\end{theorem}

\begin{theorem}\label{maintm2}
Cases 1) -- 4) are realizable by some polynomials $R$.
\end{theorem}

The theorems are proved in Sections \ref{prtm1} and \ref{prtm2}, respectively.

\section{Properties of the Schur-Szeg\H{o} composition\protect\label{knownfacts}}

The following formulae concerning the composition of polynomials can be checked straighforwardly (see \cite{Ko1}). In the second of them $S$ denotes a polynomial of degree $n-1$. Notice that in their left-hand sides (respectively, in their right-hand sides) the polynomials are composed as degree $n$ (respectively degree $n-1$) ones: 

\begin{equation}\label{diffpol}
 (A*B)'=\frac{1}{n}(A'*B')~~~~~,~~~~~(xS*B)=\frac{x}{n}(S*B')
\end{equation}
The analogs of these formulae in the case of entire transcendental functions read:

\begin{equation}\label{diffEF}
(f*g)'=f'*g'~~~~~,~~~~~xf*g=x(f*g')
\end{equation}

For the composition factors of the form $K_{a_i}$ the following formula holds:

\begin{equation}\label{formulaK}
K_{a_i}=(x+1)^{n-1}(x+a_i)=\sum _{j=0}^n{n\choose j}\left( \frac{n-j}{n}a_i+\frac{j}{n}\right) x^j
\end{equation}
Its analog in the case of composition factors $\kappa _{a_i}$ looks like this:

\begin{equation}\label{formulakappa}
\kappa _{a_i}:=e^x\left( 1+\frac{x}{a_i}\right) =\sum _{j=0}^{\infty} \frac{1}{j!}\left( 1+\frac{j}{a_i}\right) x^j
\end{equation}

The numbers $(n-j)a_i/n+j/n$ (see formula (\ref{formulaK})) form an arithmetic progression. Present the numbers $1+j/a_i$ 
from formula (\ref{formulakappa}) in the form $(a_i+j)/a_i$. Hence the numerators form also an arithmetic progression. 
This implies the following result:

\begin{cor}\label{corK}
When $a_i$ is real, there is at most one sign change in the sequence of coefficients 
of a composition factor $K_{a_i}$ or $\kappa _{a_i}$.
\end{cor}

The following proposition is Proposition 1.4 in \cite{KoSh}.

\begin{prop} \label{mult}
If the degree $n$ polynomials $A$ and $B$ have roots $x_A\neq 0$ and $x_B\neq 0$ of multiplicities
$m_A$ and $m_B$ respectively, where $m_A+m_B\geq n$, then $-x_Ax_B$ is a root of $A*B$
of multiplicity $m_A+m_B-n$.
\end{prop}

The conditions $x_A\neq 0$, $x_B\neq 0$ are omitted in \cite{KoSh} which is an error.

Set $K_{\infty}:=(x+1)^{n-1}$. For a degree $n$ polynomial $P$ denote by $P^R$ its {\em reverted} polynomial $x^nP(1/x)$. The following facts are straightforward:

\begin{prop}\label{revert}
{\rm (1)} If $P=K_{a_1}*\cdots *K_{a_{n-1}}$, then $P^R=K_{1/a_1}*\cdots *K_{1/a_{n-1}}$. 

{\rm (2)} For any two polynomials $A$ and $B$ one has $(A*B)^R=A^R*B^R$.

{\rm (3)} For any degree $n$ polynomial $A$ one has $(x+1)^{n-1}*A=A-xA'/n$ and $(x+1)^{n-2}*A=A-2xA'/n+x^2A''/(n(n-1))$.
\end{prop}

\section{Proofs of Propositions \protect\ref{necessarypol} and \protect\ref{necessaryex}\protect\label{proofsprops}}

$1^0$. We begin the proof with the following 
\vspace{2mm}

{\bf Observation.} (1) The polynomial $P$ has a $k$-fold root at $0$ ($k>0$) if and only if there are $k$ composition factors $K_{a_i}$ such that  
in $K_{a_i}$ the coefficient of $x^i$ is $0$. According to formula (\ref{formulaK}), one must have $a_i=b_i$, $i=1,\ldots ,k-1$, and 
$a_k=0$ (after a suitable permutation of the indices if necessary). 

(2) In the same way, the polynomial $R$ has a $k$-fold root at $0$ ($k>0$) if and only if there are $k$ composition factors 
$\kappa _{a_i}$ with $a_i=0,-1,\ldots ,-(k-1)$, see formula (\ref{formulakappa}).
\vspace{2mm}

$2^0$. For $k=0$, Remark~6 in \cite{Ko3} states that in the case of a polynomial $P$ there are at least $m$ different negative among the numbers $a_i$. In the case of a polynomial $R$ the same statement is contained in Corollary~2 in \cite{Ko7}.

$3^0$. Suppose that $k>0$. Set $P:=x(x+1)S$, where deg$S=n-2$. Then one of the numbers $a_i$ defined after the polynomial $P$  equals $0$. Consider the presentation

$$x(x+1)S=x(x+1)^{n-1}*(x+1)^{n-1}(x+a_2)*\cdots *(x+1)^{n-1}(x+a_{n-1})~.$$
Using both formulae (\ref{diffpol}) and formula (\ref{identity}) we present the right-hand side in the form

$$x[(x+1)^{n-2}(x+d_2)*\cdots *(x+1)^{n-2}(x+d_{n-1})]~~,~~{\rm where}~~d_i=\frac{(n-1)a_i+1}{n}~.$$
(The polynomials are composed as degree $n-1$ ones.) Hence the numbers $d_i$ are the numbers $a_i$ computed for the degree $n-1$ polynomial $(x+1)S$. One has $a_i=(nd_i-1)/(n-1)$ which implies that if $d_i<0$, then $a_i<0$.

Recall that the numbers $b_i$ depend on $i$ and $n$. In this proof we denote them further by $b_{i,n}$ because we need to compare them for different values of $n$. 

When passing from the numbers $a_i$ to the numbers $d_i$, the number $0$ corresponding to the factor $x$ in the last displayed formula is lost, and these of the numbers $a_i$ which equal $b_{j,n}$ change as follows: 
$b_{j,n}\mapsto b_{j-1,n-1}$ (to be checked directly). And in the same way, if a number $a_i$ is different from $b_{j,n}$ for all $j$, 
then $d_i$ is different from all $b_{j-1,n-1}$. Thus the rest of part (1) of Proposition~\ref{necessarypol} follows by finite induction on $k$. 
Part (1) of Proposition~\ref{necessaryex} is proved by analogy.

$4^0$. Prove part (2) of Proposition~\ref{necessarypol}. The composition $U$ of all composition factors $K_{a_i}$, 
where $a_i$ is either complex or negative (we denote their quantity by $\nu$), is a polynomial having an $(n-\nu )$-fold root at $(-1)$. 
This follows from Proposition~\ref{mult} applied $\nu -1$ times. Show that when composing $U$ consecutively with the composition 
factors $K_{a_i}$ with $0<a_i<1$, the number of negative roots (counted with multiplicity) of the given polynomial does not decrease. It suffices to consider the case when all negative roots are distinct, in the general case the result follows by continuity.

Indeed, one has $V:=(x+1)^{n-1}(x+a_i)*U=a_iU+(1-a_i)xU'$ (use formulae (\ref{diffpol}) and (\ref{identity}) and part (3) of Proposition~\ref{revert}). The signs of the polynomials 
$U$ and $V$ are the same at the negative roots of $U'$ and at $0$. These signs alternate. Hence $V$ has at most one negative root less than $U$. 

Notice that sgn$V=-$sgn$U'=$sgn$((1-a_i)xU')$ at the smallest real root $\alpha$ of $U$ (which is negative). Hence there is a root of $V$ between $\alpha$ and the smallest real root of $U'$ which is $>\alpha$, i.e. $V$ has at least as many negative roots as $U$.

$5^0$. To prove that composition with $a_i>1$ also does not decrease the number of negative roots of $U$, one can consider instead of $U$ and $K_{a_i}$ the reverted polynomials $x^nU(1/x)$ and $K_{1/a_i}$ using Proposition~\ref{revert}. One can skip the composition factors with $a_i=1$ due to (\ref{identity}). 

For $a_i=0$ one has $(x+1)^{n-1}x*U=xU'$. It is easy to show that only the first such composition can decrease by $1$ the number of negative roots, while the subsequent ones preserve this number. Part (2) of Proposition~\ref{necessarypol} is proved.

$6^0$. To prove part (2) of Proposition~\ref{necessaryex} (by analogy with $4^0$ -- $5^0$) one can use the following formula (derived from formulae (\ref{diffEF})):

$$ V:=e^x\left( 1+\frac{x}{a_i}\right)*e^xU=e^x\left( \left( 1+\frac{x}{a_i}\right) U+\frac{x}{a_i}U'\right) ~.$$
We prove the statement in the case when all negative roots of $U$ are distinct. In the presence of multiple roots the proof follows by 
continuity. 

The polynomial $xU'$ changes sign at the consecutive roots of $U$. Hence there is a root of $V$ between any two consecutive negative roots of $U$. The signs of $V$ are different at $0$ and at the smallest in absolute value root of $U$. They are also different at $-\infty$ and at the largest in absolute value root of $U$ (the details are left for the reader). Hence $V$ has one negative root more than $U$.~~~~~$\Box$

\section{Proof of Theorem~\ref{maintm1}\protect\label{prtm1}}

We prove first a proposition from which the theorem is deduced below.

\begin{prop}\label{xb_jprop}
The following composition (with $l\geq 1$ composition factors $(x+1)^{n-1}x$ and with $l+\mu \leq n$) 
\begin{equation}\label{comp01}
U:= (x+1)^{n-1}x*\cdots *(x+1)^{n-1}x*(x+1)^{n-1}(x+b_1)*\cdots *(x+1)^{n-1}(x+b_{\mu})~,
\end{equation}
is a polynomial with a $(\mu +1)$-fold root at $0$, with an $(n-\mu -l)$-fold root at $(-1)$ and with $l-1$ simple roots belonging to the interval $(-1,0)$. 

\end{prop}

{\em Proof:}\\ 

Set $T:=(x+1)^{n-1}x*(x+1)^{n-1}(x+b_1)*\cdots *(x+1)^{n-1}(x+b_{\mu})$. It follows from Proposition~\ref{mult} and from the Observation from 
$1^0$ of the proof of Propositions~\ref{necessarypol} and \ref{necessaryex} (see Section~\ref{proofsprops}) that $T$ is a polynomial of degree $n$, with a $(\mu +1)$-fold root at $0$ and with an $(n-\mu -1)$-fold one at $(-1)$. 

Denote by $T_k$ the result of composing $k$ times $T$ with $(x+1)^{n-1}x$. Apply the second of formulae (\ref{diffpol}) and then formula (\ref{identity}):

$$ (x+1)^{n-1}x*T_k=x((x+1)^{n-1}*T'_k)=xT'_k~.$$ 
If $T_k$ has a $(\mu +1)$-fold root at $0$, then this is the case of $xT'_k$ as well. The multiplicity of $(-1)$ decreases by $1$. If $-1<\zeta _1<\cdots <\zeta _k<0$ denote the roots of $T_k$ in the interval $(0,1)$, then in each of the intervals $(-1,\zeta _1)$, $(\zeta _1,\zeta _2)$, $\ldots$, $(\zeta _k,0)$ there is exactly one root of $xT'_k$. These roots are simple (because deg$(xT'_k)=n$) and the proposition is thus proved by finite induction on $k$.~~~~~$\Box$\\ 


{\em Proof of the theorem:}\\ 

{\em Case 1)}. 
\vspace{2mm}

$1^0$. Use the proposition with $l=q+q_1+q_{\mathbb{C}}$, $\mu =k-1+m+k_1+k_{\mathbb{C}}$. Hence the polynomial $U$ (see (\ref{comp01})) 
has $q-1+q_1+q_{\mathbb{C}}$ distinct roots belonging to $(-1,0)$, a $(k+m+k_1+k_{\mathbb{C}})$-fold root at $0$ and a simple root at $(-1)$. Perturb the composition factors $(x+1)^{n-1}x$ as follows:

-- $q$ of them do not change;

-- $q_1$ of them are replaced by composition factors $(x+1)^{n-1}(x+\varepsilon g_j)$, where $g_j$ are distinct positive numbers;

-- $q_{\mathbb{C}}$ of them are replaced by factors $(x+1)^{n-1}(x+\varepsilon h_j)$, where the numbers $h_j$ form $q_{\mathbb{C}}/2$ distinct conjugate couples.

Hence for $\varepsilon >0$ small enough the roots of the polynomial $U$ which belong to $(-1,0)$ are perturbed and its other roots do not change. The perturbed roots remain negative and distinct. 

$2^0$. Change the $\tilde{m}:=m+k_1+k_{\mathbb{C}}$ of the numbers $b_i$ with largest absolute values to $b_i+\lambda _i$, where $\lambda _i$ are small real parameters. Before the change the polynomial $U$ was of the form $x^{k+\tilde{m}}U_1$, $U_1(0)\neq 0$. After the change it becomes 

$$V:=U+x^k(w_1\lambda _1x^{\tilde{m}-1}+w_2\lambda _2x^{\tilde{m}-2}+\cdots +w_{\tilde{m}}\lambda _{\tilde{m}}+P)~,$$
where $w_i$ are non-zero real numbers and $P$ is a polynomial in $x$, $\lambda _1$, $\ldots$, $\lambda _{\tilde{m}}$ of total degree $\tilde{m}$ and only with monomials whose total degree w.r.t. the variables $\lambda _i$ is $\geq 2$. 

The polynomial $V/x^k$ is a versal deformation of the germ of a function $U_1$ at $0$ which has a root at $0$ of multiplicity $\tilde{m}$. Hence one can choose the values of the parameters $\lambda _i$ such that this $\tilde{m}$-fold zero splits into $m$ positive, $k_1$ negative roots and $k_{\mathbb{C}}/2$ complex conjugate couples of roots. This proves the theorem in Case 1).
\vspace{2mm}

{\em Case 2)}. 
\vspace{2mm}

Use Case 1) of the theorem with $k=q=1$. Exactly one of the composition factors equals $(x+1)^{n-1}x$. Perturb it into $(x+1)^{n-1}(x-\varepsilon )$ ($\varepsilon >0$). This perturbs the existing roots of the polynomial and its root at $0$ becomes a simple positive root (its sign can be deduced from the sign of the constant term of the polynomial which is the sign of the product of all numbers $a_i$). The existing roots were simple, therefore they remain simple after the perturbation. The numbers of negative and positive  (perturbed existing) roots remain the same. Thus one obtains instead of the figure describing Case 1) the one describing Case 2) with 
$m$ replaced by $m+1$. Hence the possibility to have $m=0$ in Case 2) has to be considered separately.

In Case 2) with $m=0$ one can again use Case 1) with $k=1$, but this time one perturbs the composition factor $(x+1)^{n-1}x$ into $(x+1)^{n-1}(x+\varepsilon )$, $\varepsilon >0$. Thus one obtains Case 2) with $q_1$ replaced by $q_1+1$. So one has to consider separately the possibility $m=q_1=0$.

For $m=q_1=0$ apply Proposition~\ref{xb_jprop} with $l=q_{\mathbb{C}}+1$, $\mu =k_1+k_{\mathbb{C}}-1$. Hence $U$ has $q_{\mathbb{C}}$ negative simple roots and a $(\mu +1)$-fold root at $0$. Perturb the $k_{\mathbb{C}}$ composition factors with largest $|b_j|$ into $(x+1)^{n-1}(x+b_j+\lambda _j)$ . (For $k_1=0$ one perturbs all $k_{\mathbb{C}}-1$ of them as indicated and one factor $(x+1)^{n-1}x$ into $(x+1)^{n-1}(x+\lambda _0)$.) The perturbation can be carried out so that the root of $U$ at $0$ split into $k_{\mathbb{C}}/2$ distinct complex conjugate couples and a $k_1$-fold root at $0$ (by analogy with Case 1), see $2^0$). 

When $k_1>0$, the $k_1-1$ remaining factors $(x+1)^{n-1}(x+b_j)$ and one factor $(x+1)^{n-1}x$ are perturbed so that $U$ have $k_1$ negative roots close to $0$ (the previously existing $q_{\mathbb{C}}$ negative roots remain such). Its $k_{\mathbb{C}}$ complex couples remain such. Finally, perturb the remaining factors $(x+1)^{n-1}x$ into $(x+1)^{n-1}(x+h_i)$, where the 
numbers $h_i$ form $q_{\mathbb{C}}/2$ conjugate couples. 
\vspace{2mm}

{\em Case 3)} 
\vspace{2mm}

Apply Proposition~\ref{xb_jprop} with $l=q_1+q+s$, $\mu =k+r+m$. Notice that for $\delta =0$ one has $l+\mu <n$. The polynomial $U$ has a $(\mu +1)$-fold root at $0$ and $l-1$ negative roots. Perturb the $r+m$ composition factors with largest $|b_j|$ so that the root of $U$ at $0$ split into $r/2$ complex conjugate couples, a $k$-fold root at $0$ and $m$ positive roots. Then perturb the factors $(x+1)^{n-1}x$ into $(x+1)^{n-1}(x+\eta _i)$ as follows:

-- $s$ of the numbers $\eta _i$ form distinct complex conjugate couples;

-- $q_1$ of them are negative and distinct;

-- $q$ of them are $0$.

The last perturbation does not change the number of positive, negative, zero and complex roots of $U$. For $\delta =0$ the case is completely solved.

Denote by $S$ the sector $\{ u+iv\in \mathbb{C}|0<2u<v\}$. When $\delta >0$, we need the following lemma (proved after the proof of the theorem).

\begin{lm}\label{lmconstructclose}
{\rm (1)} For all $\varepsilon \in S$ sufficiently close to $0$ the polynomial 

$$V:=(x+1)^{n-1}(x+1+\varepsilon )*(x+1)^{n-1}(x+1+\bar{\varepsilon})$$ 
has two complex conjugate roots close to $(-1)$ and an $(n-2)$-fold root at $(-1)$.

{\rm (2)} Suppose that the polynomial $U$ has $s^*$ negative and $t^*$ positive simple roots, $r^*/2$ distinct conjugate couples and an $(n-s^*-t^*-r^*)$-fold root at $(-1)$, $n-s^*-t^*-r^*\geq 2$. 
Then one can choose $\varepsilon \in S$ so close to $0$ that the polynomial $U*V$ have $s^*$ negative (different from $-1$) and $t^*$ positive simple roots, $r^*/2+1$ distinct complex conjugate couples and an $(n-s^*-t^*-r^*-2)$-fold root at $(-1)$. 

{\rm (3)} The multiplicity of $0$ as a root of $U$ and $U*V$ is the same for $\varepsilon$ small enough.
\end{lm}

To complete the proof in Case 3) one sets $s^*:=q-1+q_1+s$, $t^*:=m$, $r^*:=r$ and then applies the lemma $\delta /2$ times. 
\vspace{2mm}

{\em Case 4)}
\vspace{2mm}

Use Case 3) with $k=q=1$. Perturb the factor $(x+1)^{n-1}x$ into $(x+1)^{n-1}(x-\zeta )$, $\zeta >0$. This changes $m$ to $m+1$. Thus Case 
4) is deduced from Case 3) except for $m=0$. For $m=0$ use again Case 3) with $k=q=1$ changing this time $(x+1)^{n-1}x$ into 
$(x+1)^{n-1}(x+\zeta )$. This changes $q_1$ into $q_1+1$ and there remains to consider only the possibility $m=q_1=0$.  

In this particular case $r$ must be even. Hence such are $\delta$ and $s$ as well. Apply Proposition~\ref{xb_jprop} with $l=s+1$, $\mu =r-1$. Hence the polynomial $U$ has $s$ negative distinct roots and an $r$-fold root at $0$. Perturb one factor $(x+1)^{n-1}x$ and the factors $(x+1)^{n-1}(x+b_j)$ to make the $r$-fold root of $U$ split into $r/2$ distinct conjugate couples. Hence the factor $(x+1)^{n-1}x$ becomes $(x+1)^{n-1}(x-\varepsilon )$ with $\varepsilon >0$ (this follows from $r$ being even). After this perturb the remaining $s$ factors $(x+1)^{n-1}x$ into $(x+1)^{n-1}(x+\zeta _i)$, where the numbers $\zeta _i$ form $s/2$ distinct conjugate couples. This finishes the construction for $\delta =0$. For $\delta >0$ 
one has to apply $\delta /2$ times Lemma~\ref{lmconstructclose} with $r^*=r$, $t^*=0$ and $s^*=q_1+s$.~~~~~$\Box$\\ 

{\em Proof of Lemma~\ref{lmconstructclose}:}\\ 

Using equality (\ref{identity}) and part (3) of Proposition~\ref{revert} one finds that 

$$V=(x+1)^{n-2}((x+1)^2+(\varepsilon +\bar{\varepsilon}+\varepsilon\bar{\varepsilon}/n)(x+1)+(n-1)\varepsilon\bar{\varepsilon}/n)~.$$
Set $\varepsilon =u+iv$. The discriminant $\Delta$ of the quadratic factor (considered as a polynomial in $x+1$) equals 

$$(4/n)(u^2-(n-1)v^2)+o(u^2+v^2)~.$$ 
When $\varepsilon \in S$ 
is close to $0$, one has $\Delta <0$. The coefficient of $(x+1)$ and the constant term of the quadratic factor tend to $0$ as $\varepsilon \rightarrow 0$. Hence its roots also tend to $0$. This proves part (1) of the lemma.

Set $p:=n-s^*-t^*-r^*$. Present $U$ as a polynomial in $x+1$: 

$$U=U_0(x+1)^{p}+\cdots +U_{n-p}(x+1)^n~,~U_0\neq 0~.$$ 
One has 

$$U*V=U_0(x+1)^{p}*V+\cdots +U_{n-p}(x+1)^n*V~.$$
By Proposition~\ref{mult} all terms have a root at $(-1)$ of multiplicity at least $p-2$. As $V$ is a perturbation of $(x+1)^n$, for $\varepsilon$ small enough the polynomial $U*V$ has $r^*/2$ conjugate couples, $s^*$ negative and $t^*$ positive roots close to the ones of $U$. We show that when $\varepsilon \in S$ is small enough, then the first of the terms to the right has two complex conjugate roots close to $-1$. One has 

$$(x+1)^{p}*V=(x+1)^{p-2}((1+o(1))(x+1)^2+(2pu/n+o(|u|+|v|))(x+1)+$$
$$((u^2+v^2)p(p-1)/n^2))~.$$
The discriminant of the quadratic factor equals 

$$(4/n^2)(p^2u^2-p(p-1)(u^2+v^2))+o(u^2+v^2)~.$$ 
For $\varepsilon \in S$ small enough it is $<0$. The coefficient of $x+1$ and the constant term of the quadratic factor tend to $0$ as $\varepsilon \rightarrow 0$ while the coefficient of $(x+1)^2$ remains close to $1$, so its two roots also tend to $-1$.  

Set $\varepsilon =\tau \eta$, where $\tau =|\varepsilon |$. Set $x+1=\tau y$. Set $L:=(V/(x+1)^{n-2})|_{x=\tau y-1}$ and $B:=((x+1)^{p}*V/(x+1)^{p-2})|_{x=\tau y-1}$. As $\tau \rightarrow 0$, the two complex roots of the polynomial $L$ (respectively $B$) are of the form $\alpha _i\tau +o(\tau )$, (respectively $\beta _i\tau +o(\tau )$), $\alpha _i\neq 0\neq \beta _i$, $i=1,2$. The numbers $\alpha _i$ and $\beta _i$ are roots respectively of the polynomials

$$y^2+(\eta +\bar{\eta})y+((n-1)/n)\eta \bar{\eta}~~~~~~~~\hspace{5mm}{\rm and}$$

$$y^2+(p(\eta +\bar{\eta})/n)y+(p(p-1)/n^2)\eta \bar{\eta}~.$$

Consider two circles of radius $\tau \min (|\beta _1|, |\beta _2|)/2$ centered at the roots of $B$. All terms 
$(U_j(x+1)^{p+j}*V)|_{x=\tau y-1}$, $j=1,\ldots ,n-p$ when restricted to these circles have their module tending to $0$ (as $\tau \rightarrow 0$) faster than the module of the term $U_0B$. (This is due to their higher power of $(x+1)$, i.e. of $\tau y$.) By Rouch\'e's theorem inside each of the circles there is exactly one root of the polynomial $U*V$. This proves part (2) of the lemma. 

Part (3) is evident -- the multiplicity of $0$ as a root of $U$ and $U*V$ is defined by the number of the first consecutive coefficients of these polynomials which are $0$. As $V$ is a perturbation of $(x+1)^n$, these numbers are the same.~~~~~$\Box$

\section{Proof of Theorem~\protect\ref{maintm2}\protect\label{prtm2}}

The theorem is proved with the help of the following proposition:

\begin{prop}\label{xb_jpropexp}
The following composition (with $l\geq 1$ composition factors $xe^x$) 
\begin{equation}\label{comp01exp}
U:= e^xx*\cdots *e^xx*e^x(x-1)*\cdots *e^x(x-\mu )~,
\end{equation}
is of the form $e^xY$, where $Y$ is a degree $l+\mu$ polynomial with a $(\mu +1)$-fold root at $0$ and with $l-1$ simple negative roots. 
\end{prop}

{\em Proof:}\\ 

With the help of formulae (\ref{diffEF}) and using finite induction on $\mu$ one shows that the composition of the last $\mu +1$ composition factors is exactly $e^xx^{\mu +1}$. 

Suppose that the proposition is true for $l=l_0$. Then for $l=l_0+1$ one has $U=e^xx*e^xY=e^xx(Y+Y')$. The sign of $x(Y+Y')$ changes alternatively at the consecutive negative roots of $Y'$. Hence there is a root of $U$ between any two negative roots of $Y'$. Denote the latter roots by 
$\beta _1<\cdots <\beta _{l_0-1}$ and by $\alpha$ the greatest (i.e. smallest in absolute value) negative root of $Y$.

One has sgn$U(-\infty )=-$sgn$Y(-\infty )$, sgn$U(\beta _1)=-$sgn$Y(\beta _1)$ and sgn$Y(-\infty )=-$sgn$Y(\beta _1)$. Hence there is a root of $U$ in $(-\infty ,\beta _1)$. In the same way, sgn$U(\beta _{l_0-1})=-$sgn$Y(\beta _{l_0-1})$, sgn$U(\alpha )=-$sgn$Y'(\alpha )=$sgn$Y(\beta _{l_0-1})$. Hence $U$ has a root in $(\beta _{l_0-1},\alpha )$ as well.

Thus the product $x(Y+Y')$ has $l_0$ distinct negative roots and a $(\mu +1)$-fold root at $0$. As deg$(Y+Y')=l_0+\mu$, all negative roots are simple.~~~~~$\Box$ \\

{\em Proof of the theorem:}\\

{\em Case 1)}. 
\vspace{2mm}

Use the proposition with $l=q+q_1+q_{\mathbb{C}}$, $\mu =k-1+m+k_1+k_{\mathbb{C}}$. Hence the polynomial $Y$ 
has $q-1+q_1+q_{\mathbb{C}}$ distinct negative roots and a $(k+m+k_1+k_{\mathbb{C}})$-fold root at $0$. Perturb the composition factors $e^xx$ as follows:

-- $q$ of them do not change;

-- $q_1$ of them are replaced by composition factors $e^x(x+\varepsilon g_j)$, where $g_j$ are distinct positive numbers;

-- $q_{\mathbb{C}}$ of them are replaced by factors $e^x(x+\varepsilon h_j)$, where the numbers $h_j$ form $q_{\mathbb{C}}/2$ distinct conjugate couples.

Hence for $\varepsilon >0$ small enough the negative roots of the polynomial $Y$ are perturbed and its other roots do not change. The perturbed roots remain negative and distinct. 

Change the $\tilde{m}:=m+k_1+k_{\mathbb{C}}$ of the composition factors $e^x(x-j)$ with largest absolute values of $j$ to $e^x(x-j+\lambda _j)$, where $\lambda _j$ are small real parameters chosen such that the root of $Y$ at $0$ split into a $k$-fold root at $0$, $k_{\mathbb{C}}/2$ complex conjugate couples and $k_1$ negative roots. The proof of Case 1) is finished by complete analogy with the proof of this case in Theorem~\ref{maintm1}. 

The rest of the proof of Theorem~\ref{maintm2} is done also by analogy with the rest of the proof of Theorem~\ref{maintm1} (modulo some technical details) -- the role of the composition factor $(x+1)^{n-1}x$ in the latter is played by $e^xx$, the one of $(x+1)^{n-1}(x+b_j)$ is played by $e^x(x-j)$. 
\vspace{2mm}

{\em Case 2)} 
\vspace{2mm}

Use Case 1) of the theorem with $k=q=1$. Exactly one of the factors equals $e^xx$. Perturb it into $e^x(x-\varepsilon )$ ($\varepsilon >0$). The root at $0$ of the polynomial becomes a simple positive root. The existing roots remain simple after the perturbation. The numbers of negative and positive  (perturbed existing) roots remain the same. Thus one obtains (instead of Case 1)) Case 2) with 
$m$ replaced by $m+1$. The possibility to have $m=0$ in Case 2) has to be considered separately.

For $m=0$ one can again use Case 1) with $k=1$, this time perturbing the factor $e^xx$ into $e^x(x+\varepsilon )$, $\varepsilon >0$. Thus one obtains Case 2) with $q_1$ replaced by $q_1+1$. There remains to consider the possibility $m=q_1=0$.

For $m=q_1=0$ apply Proposition~\ref{xb_jpropexp} with $l=q_{\mathbb{C}}+1$, $\mu =k_1+k_{\mathbb{C}}-1$. Hence $Y$ has $q_{\mathbb{C}}$ negative roots and a $(\mu +1)$-fold root at $0$. Perturb the $k_{\mathbb{C}}$ composition factors $e^x(x-j)$ with largest $j$ into $e^x(x-j+\lambda _j)$. (When $k_1=0$ one perturbs all $k_{\mathbb{C}}-1$ factors $e^x(x-j)$ as indicated and one factor $e^xx$ into $e^x(x+\lambda _0)$.) The root of $Y$ at $0$ splits into $k_{\mathbb{C}}/2$ distinct complex conjugate couples and a $k_1$-fold root at $0$. 

When $k_1>0$, the $k_1-1$ remaining factors $e^x(x-j)$ and one factor $e^xx$ are perturbed so that $Y$ have $k_1$ negative roots close to $0$ (the previously existing $q_{\mathbb{C}}$ negative roots remain such). Its $k_{\mathbb{C}}$ complex couples remain such. Finally, perturb the remaining factors $e^xx$ into $e^x(x+h_i)$, where the numbers $h_i$ form $q_{\mathbb{C}}/2$ conjugate couples.
\vspace{2mm}

{\em Case 3)} 
\vspace{2mm}

Apply Proposition~\ref{xb_jpropexp} with $l=q_1+q+s$, $\mu =k-1+r+m$. (For $\delta =0$ one has $l+\mu =n$, otherwise $l+\mu <n$.) The polynomial $Y$ has a $(\mu +1)$-fold root at $0$ and $l-1$ negative roots. Perturb the $r+m$ composition factors $e^x(x-j)$ with largest $j$ so that the root of $Y$ at $0$ split into $r/2$ complex conjugate couples, a $k$-fold root at $0$ and $m$ positive roots. Then perturb the factors $e^xx$ into $e^x(x+\eta _i)$ as follows:

-- $s$ of the numbers $\eta _i$ form distinct complex conjugate couples;

-- $q_1$ of them are negative and distinct;

-- $q$ of them are $0$.

The last perturbation does not change the number of positive, negative, zero and complex roots of $Y$. For $\delta =0$ this finishes the proof of Case 3). 

Denote by $S$ some sector centered at $0$ and avoiding (except $0$) the real axis. For $\delta >0$ we need the following lemma (see its proof after the proof of the theorem).

\begin{lm}\label{LLLL}
{\rm (1)} One has $e^x(1+\varepsilon x)*e^x(1+\bar{\varepsilon}x)=e^xV$, where 

$$V=1+(\varepsilon +\bar{\varepsilon}+\varepsilon\bar{\varepsilon})x+\varepsilon\bar{\varepsilon}x^2~.$$

{\rm (2)} Suppose that the polynomial $Y$ is of degree $p$. 
Then the function $(e^xY)*(e^xV)$ is of the form $e^xY_1$, where $Y_1$ is a degree $p+2$ polynomial. 

{\rm (3)} Suppose that the degree $p$ polynomial $Y$ has $s^*$ negative and $t^*$ positive simple roots, $r^*/2$ distinct conjugate couples and a $(p-s^*-t^*-r^*)$-fold root at $0$. One can choose $\varepsilon \in S$ so close to $0$ that the polynomial $Y_1$ have $s^*$ negative and $t^*$ positive simple roots, $r^*/2+1$ distinct conjugate couples and a $(p-s^*-t^*-r^*)$-fold root at $0$. 
\end{lm}

Set $r^*:=r$, $t^*:=m$ and $s^*:=q-1+q_1+s$. Applying the lemma $\delta /2$ times one obtains the proof of the theorem in Case 3).
\vspace{2mm}

{\em Case 4)} 
\vspace{2mm}

Use Case 3) with $k=q=1$. Perturb the factor $e^xx$ into $e^x(x-\zeta )$, $\zeta >0$. This changes $m$ to $m+1$. Thus Case 
4) is deduced from Case 3) except for $m=0$. For $m=0$ use again Case 3) with $k=q=1$ changing this time $e^xx$ into 
$e^x(x+\zeta )$. This changes $q_1$ into $q_1+1$ and there remains to consider only the possibility $m=q_1=0$.  

In this particular case $r$, $\delta$ and $s$ are even. Apply Proposition~\ref{xb_jpropexp} with $l=s+1$, $\mu =r-1$. Hence the polynomial $Y$ has $s$ negative roots and an $r$-fold root at $0$. Perturb one factor $e^xx$ and the factors $e^x(x-j)$ to make the $r$-fold root of $Y$ split into $r/2$ conjugate couples. Hence the factor $e^xx$ becomes $e^x(x-\varepsilon )$ with $\varepsilon >0$ (this follows from $r$ being even). After this perturb the remaining $s$ factors $e^xx$ into $e^x(x+\delta _i)$, where the numbers $\delta _i$ form $s/2$ distinct conjugate couples. This finishes the construction for $\delta =0$. 

For $\delta >0$ one has to apply $\delta /2$ times Lemma~\ref{LLLL}.~~~~~$\Box$\\ 

{\em Proof of Lemma~\ref{LLLL}:}\\ 

Parts (1) and (2) of the lemma follow from the second of formulae (\ref{diffEF}) with $f=e^xx$ or $f=e^xx^2$ and $g=e^xY$. Prove part (3). The polynomial $Y_1$ is a perturbation of the polynomial $Y$, therefore for $\varepsilon$ close to $0$ it has $p$ roots close to the respective roots of $Y$ and two roots (called {\em distant}) whose moduli tend to $\infty$ as $\varepsilon \rightarrow 0$. 

The polynomials $Y$ and $Y_1$ have the same multiplicity of the root at $0$. Indeed, for $\varepsilon$ nonreal all coefficients of the function $e^xV$ are nonzero and this multiplicity is defined by the number of first consecutive coefficients of $e^xY$ which are $0$. 

As both $Y$ and $Y_1$ are real polynomials, $Y_1$ has the same number of distinct negative and distinct positive roots and the same number of distinct complex conjugate couples as $Y$ (excluding the two distant roots). 

Suppose that $Y$ is monic (this is not restrictive). For $p^*\in \mathbb{N}$ one has 

$$e^xV*e^xx^{p^*}=x^{p^*}(\varepsilon\bar{\varepsilon}x^2+(\varepsilon +\bar{\varepsilon}+(2p^*+1)\varepsilon\bar{\varepsilon})x+1+p^*(\varepsilon +\bar{\varepsilon})+{p^*}^2\varepsilon\bar{\varepsilon})~.$$
Set $W:=\varepsilon\bar{\varepsilon}x^2+(\varepsilon +\bar{\varepsilon})x+1$. Hence the function $e^xV*e^xY$ is of the form 
$e^xx^p(W+T)$, where $T$ is a Laurent series in $x$ whose coefficients are polynomials in $\varepsilon$ and $\bar{\varepsilon}$. It contains only monomials $x^{\alpha}\varepsilon ^{\beta}\bar{\varepsilon}^{\gamma}$ with $\alpha -\beta -\gamma <0$. 

The roots of $W$ are $1/\varepsilon$ and $1/\bar{\varepsilon}$. Consider two circles $C_1$ and $C_2$ centered at them and of radius $1$. When $\varepsilon \in S$ is small, the values of $|T|$ at each point of each of the two circles are much smaller than the respective values of $W$. By Rouch\'e's theorem each of the circles contains exactly one root of $x^p(W+T)$.~~~~~$\Box$


\begin{thebibliography}{Dillo 83}
\bibitem{AlKo} S. Alkhatib and V.P. Kostov,
The Schur-Szeg\"o composition of real polynomials of degree $2$,
Rev. Mat. Complutense 21 (2008) no. 1, 191--206.






\bibitem{Ko1} V. P. Kostov, The Schur-Szeg\"o composition for 
hyperbolic polynomials, C.R.A.S. S\'er. I 345/9 (2007), 483-488, 
doi:10.1016/j.crma.2007.10.003.

\bibitem{Ko2} V. P. Kostov, Eigenvectors in the context of the 
Schur-Szeg\"o composition of polynomials, Math. Balkanica 
22(2008) Fasc. 1-2, 155--173.

\bibitem{Ko3} V. P. Kostov, Teorema realizatsii v kontekste kompozitsii 
Shura-Sege, Funkcional'nyy Analiz i ego Prilozheniya 43 (2009) no. 2, 79-83. 
(A realization theorem in the context of the 
Schur-Szeg\"o composition, Funct. Anal. Appl. 43 (2009) no. 2, 147-150.)

\bibitem{Ko4} V. P. Kostov, A mapping connected with the Schur-Szeg\H{o} 
composition, C.R.A.S. S\'er. I 347 (2009) 1355-1350.

\bibitem{Ko5} V.P. Kostov, The Schur-Szeg\"o composition for real polynomials, 
C.R.A.S. S\'er. I, 346 (2008), 271-276.

\bibitem{Ko6} V.P. Kostov, A realization theorem about the Schur-Szeg\"o 
composition for entire functions, Comptes Rendus Acad. Sci. Bulgare 
62, No. 1 (2009), 17-22. 

\bibitem{Ko7} V. P. Kostov, A mapping defined by the Schur-Szeg\H{o} 
composition, Comptes Rendus Acad. Sci. Bulgare Vol. 63 (2010) No. 7 943-952.

\bibitem{KoSh} V. P. Kostov and B. Z. Shapiro, 
On the Schur-Szeg\"o composition of polynomials, C.R.A.S. S\'er. I 
343 (2006) 81--86. 

\bibitem{KoShMa} V. P. Kostov, B. Z. Shapiro and A. Martinez-Finkelstein, 
Narayana numbers and Schur-Szeg\"o composition,  
J. Approx. Theory, 161 (2) (2009) 464-476.



\bibitem{Pr} V. Prasolov, Polynomials, Translated from the 2001 
Russian second edition by Dimitry Leites. Algorithms and Computation 
in Mathematics, 11. Springer-Verlag, Berlin, 2004.

\bibitem{RS} Q. I. Rahman and G. Schmeisser, Analytic Theory of Polynomials, 
London Math. Soc. Monogr. (N.S.), vol. 26, Oxford Univ. Press, New 
York, NY, 2002. 
\end{thebibliography}
\end{document}